\documentclass[11pt,reqno]{amsart}   
\usepackage{amssymb}
\usepackage{amsthm}
\usepackage{amsmath}
\usepackage{amscd}
\usepackage{comment}
\usepackage{euscript}

\newtheorem{thm}{Theorem}[section]
\newtheorem{prop}[thm]{Proposition}

\newtheorem{lemma}[thm]{Lemma}

\theoremstyle{definition}
\newtheorem{defi}[thm]{Definition}
\newtheorem{notation}[thm]{Notation}

\theoremstyle{remark}
\newtheorem{rk}[thm]{Remark}
\newtheorem{obs}[thm]{Observation}
\newtheorem{ques}[thm]{Question}

\newcount\theTime
\newcount\theHour
\newcount\theMinute
\newcount\theMinuteTens
\newcount\theScratch
\theTime=\number\time
\theHour=\theTime
\divide\theHour by 60
\theScratch=\theHour
\multiply\theScratch by 60
\theMinute=\theTime
\advance\theMinute by -\theScratch
\theMinuteTens=\theMinute
\divide\theMinuteTens by 10
\theScratch=\theMinuteTens
\multiply\theScratch by 10
\advance\theMinute by -\theScratch

\def\today{{\number\day\space
 \ifcase\month\or
  January\or February\or March\or April\or May\or June\or
  July\or August\or September\or October\or November\or December\fi
 \space\number\year}}

\newcommand{\NN}{{\mathbb{N}}}
\newcommand{\ZZ}{{\mathbb{Z}}}

\newcommand{\CC}{{\mathbb{C}}}

\newcommand{\FF}{{\mathbb{F}}}
\newcommand{\cstar}{\mbox{$C^*$}}

\newcommand\HEu{{\EuScript H}}                   
\newcommand\HEut{{\widetilde\HEu}}                   

\newcommand{\fl}{{\mathcal{L}}}
\newcommand{\fk}{{\mathcal{K}}}

\newcommand{\ha}{{\hat{a}}}

\newcommand{\tdt}{\otimes\dots\otimes}

\newcommand\sigmat{{\tilde\sigma}}
\newcommand\id{{\operatorname{id}}}

\newcommand\Aut{{\operatorname{Aut}}}
\newcommand\eps{\epsilon}

\title[Unique ergodicity]{Unique ergodicity of free shifts and some other automorphisms of C$^*$--algebras} 
\author[Abadie, Dykema]{Beatriz Abadie*, Ken Dykema$^\dag$}
\address{Dept. of Mathematics, Texas A\&M University, College Station TX 77843-3368. USA.\\
Permanent address: Centro de Matem\'aticas, Facultad de Ciencias, Igu\'a 4225, CP 11 400, Montevideo, Uruguay.}
\email{abadie@cmat.edu.uy}
\thanks{*Partially supported by Proyecto Clemente Estable 8013.
$^\dag$Partially supported by NSF grants DMS--0300336 and DMS--0600814.}

\address{Dept. of Mathematics, Texas A\&M University, College Station TX 77843-3368, USA.}
\email{Ken.Dykema@math.tamu.edu}
\date{\today}
\begin{document}



\begin{abstract}
A notion of unique ergodicity relative to the fixed--point subalgebra
is defined for automorphisms of unital \cstar--algebras.
It is proved that the free shift on any reduced amalgamated free product \cstar--algebra
is uniquely ergodic relative to its fixed--point subalgebra,
as are
autormorphisms of reduced group \cstar--algebras arising from certain
automorphisms of groups.
A generalization of Haagerup's inequality, yielding bounds on the norms of certain elements
in reduced amalgamated free product \cstar--algebras, is proved.
\end{abstract}

\maketitle


\section{Introduction}

Let $\Omega$ be a compact Hausdorff space and $T$ a homeomorphism of $\Omega$ onto itself.
In the terminology of~\cite{Ox},
(see also~\cite{KB} and~\cite{Fu}, where slightly different terminology is used),
$T$ is called
{\em uniquely ergodic} if there is a unique $T$--invariant Borel probability measure $\mu$ on $\Omega$,
(with respect to which  $T$ is then necessarily ergodic).
Oxtoby shows~\cite[5.1]{Ox} that if $T$ is uniquely ergodic, then the ergodic averages
\begin{equation}
\frac1n\sum_{k=0}^{n-1}f\circ T^k
\end{equation}
converge uniformly to the constant $\int f\,d\mu$, as $n\to\infty$.

The homeomorphisms of $\Omega$ are in 1--1 correspondence with the automorphisms of
the C$^*$--algebra $C(\Omega)$ of all continuous, complex--valued functions on $\Omega$
and the Borel probability measures on $\Omega$ are by Riesz's Theorem
in 1--1 correspondence with the states of $C(\Omega)$.
There is a natural noncommutative version of unique ergodicity.
Let $A$ be a unital C$^*$--algebra and let $\alpha$ be an automorphism of $A$.
An $\alpha$--invariant state of $A$ always exists, and can be found, for example, by taking a weak limit
of averages
\begin{equation}
\frac1n\sum_{k=0}^{n-1}\phi\circ\alpha^k
\end{equation}
of any state $\phi$.
We say $\alpha$ is {\em uniquely ergodic} if there is a unique $\alpha$--invariant state of $A$.
It is not difficult to show (based on Oxtoby's argument~\cite[5.1]{Ox}) that $\alpha$ is uniquely
ergodic if and only if for every $a\in A$ the ergodic averages
\begin{equation}\label{eq:ergodicaves}
\frac1n\sum_{k=0}^{n-1}\alpha^k(a)
\end{equation}
converge in norm to a scalar multiple of the identity as $n\to\infty$
and, in this case, the invariant state evaluated at $a$ is equal to this limit.
(A more general result is proved in Theorem~\ref{thm:ue} below.)

Our interest in these topics grew out of a question asked by David Kerr~\cite{K}:
Is the free shift on $C^*_r(\FF_\infty)$ uniquely ergodic?
A positive answer to Kerr's question follows from Haagerup's inequality~\cite{haag}.
This argument is described in section~\ref{sec:CFinf} below.

In considering more general free shift automorphisms, we were motivated to consider
a broader notion of unique ergodicity.
If $A$ is a unital C$^*$--algebra and $\alpha$ an automorphism of $A$,
consider the fixed--point subalgebra
\begin{equation}\label{eq:fps}
A^\alpha=\{a\in A\mid\alpha(a)=a\}.
\end{equation}
We say that $\alpha$ is {\em uniquely ergodic relative to its fixed--point subalgebra}
if every state of $A^\alpha$ has a unique $\alpha$--invariant state extension to $A$.
In the case when $A^\alpha$ consists only of scalar multiples of the identity element, this reduces
to the usual notion of unique ergodicity.
In section~\ref{sec:ue}, we give some alternative characterizations of unique ergodicity
relative to the fixed--point subalgebra.
It turns out to be equivalent to norm convergence of the ergodic averages~\eqref{eq:ergodicaves}
as $n\to\infty$ for all $a\in A$.
Thus, unique ergodicity relative to the fixed--point subalgebra
implies (by taking the limit of the ergodic averages)
existence of a unique
$\alpha$--invariant conditional expectation from $A$ onto $A^\alpha$.
However (see Question~\ref{qn}) we do not know whether the converse direction holds.

After seeing that the free shift on $C^*_r(\FF_\infty)$ is uniquely ergodic, it is natural to
ask whether free shifts on other reduced free product \cstar--algebras and even on
reduced amalgamated free product \cstar--algebras are uniquely ergodic
relative to their fixed point subalgebras.
We give an affirmative answer in Theorem~\ref{thm:afp}.

A technical result that we use is an extension of Haagerup's inequality to the setting of reduced
amalgamated free product \cstar--algebras.
Haagerup's inequality says that the operator norm of an element of $C^*_r(\FF_\infty)$ that
is supported on words of length $n$ is no greater than $n+1$ times the $\ell^2$--norm.
It is a fundamental inequality, and has been generalized in several different directions;
see, for example, \cite{dlH}, \cite{J}, \cite{DHR}, \cite{KS}, \cite{RX}. 
One such generalization is~\cite[3.3]{DHR}, in the context of reduced free product C$^*$--algebras
with amalgamation over the scalars,
which applies to all finite linear combinations of words of fixed block length $n$.
A strong generalization, due to Ricard and Xu~\cite{RX},
is in the context of reduced amalgamated free product
C$^*$--algebras; they prove bounds on operator norms that apply
to all matrices over all finite linear combinations of words of fixed block length $n$.
In Proposition~\ref{prop:ineq}, we prove a generalization of Haagerup's inequality in the setting
of reduced amalgamated free product C$^*$--algebras.
Our bound on the operator norm applies only to certain linear combinations of words of 
block length $n$, but our bound has a rather nice form.
In fact, as Eric Ricard kindly showed us, our Proposition~\ref{prop:ineq}
follows from the results of Ricard and Xu.
However, we nonetheless present our direct proof here, as it is slightly simpler
(for being a more specific result).

To summarize the contents:
section~\ref{sec:CFinf} contains the proof of unique ergodicity of the free shift on
$C^*_r(\FF_\infty)$;
section~\ref{sec:ue} gives alternative characterizations of unique ergodicity relative to the
fixed--point subalgebra, and contains a generalization of the
argument from the previous section to groups with property~(RD) of Jolissaint;
section~\ref{sec:constr} recalls the construction of the reduced amalgamated free product
of \cstar--algebras;
section~\ref{sec:normest} contains a generalization of Haagerup's inequality to reduced amalgamated
free product \cstar--algebras;
section~\ref{sec:freeshift} proves that free shifts are uniquely ergodic relative to their fixed--point
subalgebras.

\medskip

{\bf{Acknowledgement.}} This work was carried out while the first author was
visiting the Mathematics Department of Texas A\&M University.
She would like to thank the members of the department  for their warm hospitality.
The authors thank
Thierry Fack for a helpful comment that lead to Proposition~\ref{prop:Fack},
and David Kerr and Eric Ricard for helpful conversations.

\section{The free shift on $C^*_r(\FF_\infty)$ is uniquely ergodic}
\label{sec:CFinf}

Here, $C^*_r(\FF_\infty)$ is the reduced group \cstar--algebra of the free group on infinitely many generators
$\{g_i\}_{i\in\ZZ}$ and the free shift is the automorphism $\alpha$ of $C^*_r(\FF_\infty)$ arising
from the automorphism of the group that sends $g_i$ to $g_{i+1}$.

The \cstar--algebra $C^*_r(\FF_\infty)$ is densely spanned by the left translation operators $\lambda_h$
acting on $\ell^2(\FF_\infty)$, ($h\in\FF_\infty$).
If $h=e$ is the trivial group element, then $\lambda_h$ is the identity element $1$ and
\begin{equation}
\frac1n\sum_{k=0}^{n-1}\alpha^k(1)=1
\end{equation}
for all $n$.
If $h$ is a nontrivial element of word length $p$, then by Haagerup's inequality~\cite[1.4]{haag},
\begin{equation}\label{eq:haagineq}
\bigg\|\frac1n\sum_{k=0}^{n-1}\alpha^k(\lambda_h)\bigg\|
\le(p+1)\bigg\|\frac1n\sum_{k=0}^{n-1}\alpha^k(\lambda_h)\bigg\|_2
=\frac{p+1}{\sqrt n},
\end{equation}
where $\|\cdot\|_2$ refers to the norm of the corresponding element in $\ell^2(\FF_\infty)$.
We conclude that the averages appearing on the left--hand--side of~\eqref{eq:haagineq} tend
to zero as $n\to\infty$, and this proves that the free shift is uniquely ergodic and
that its unique invariant state is the canonical tracial state $\tau$ defined by
\begin{equation}
\tau(\lambda_h)=
\begin{cases}
1,&h=e \\
0,&h\ne e.
\end{cases}
\end{equation}

\section{Unique ergodicity relative to the fixed--point subalgebra}
\label{sec:ue}

In this section, we prove certain conditions equivalent to unique ergodicity
relative to the fixed--point subalgebra.

\begin{obs}\label{obs:Jordan}
Let $A$ be a C$^*$--algebra and let $\phi:A\to\CC$ be a self--adjoint linear functional,
namely a bounded linear functional such that $\phi(a^*)$ is the complex conjugate of $\phi(a)$.
Recall (see~\cite[3.2.5]{Ped}) that the Jordan decomposition of $\phi$ is the unique pair
$\phi_+$ and $\phi_-$ of positive linear functionals such that $\phi=\phi_+-\phi_-$ and
$\|\phi\|=\|\phi_+\|+\|\phi_-\|$.
Suppose $\alpha\in\Aut(A)$ and $\phi$ is $\alpha$--invariant.
Then $\phi=\phi\circ\alpha=\phi_+\circ\alpha-\phi_-\circ\alpha$ and
$\|\phi\|=\|\phi_+\|+\|\phi_-\|=\|\phi_+\circ\alpha\|+\|\phi_-\circ\alpha\|$.
By uniqueness, it follows that $\phi_+$ and $\phi_-$ are both $\alpha$--invariant.
\end{obs}

Recall that a {\em conditional expectation} from
a C$^*$--algebra $A$ onto a C$^*$--subalgebra $B$
is a projection $E$ of norm $1$ from $A$ onto $B$.
A classical result of Tomiyama~\cite{Tom}
is that such a projection $E$ is automatically completely positive and satisfies
the conditional expectation property.

\begin{thm}\label{thm:ue}
Let $\alpha$ be an automorphism of a unital C$^*$--algebra $A$
and let $A^\alpha$ be its fixed--point subalgebra as in~\eqref{eq:fps}.
Then the following five statements are equivalent:
\begin{itemize}
\item[(i)] Every bounded linear functional on $A^\alpha$ has a unique bounded, $\alpha$--invariant
linear extension to $A$.

\item[(ii)] Every state of $A^\alpha$ has a unique $\alpha$--invariant state extension to $A$.

\item[(iii)] The subspace $A^\alpha+\{a-\alpha(a)\mid a\in A\}$ is dense in $A$.

\item[(iv)] For all $a\in A$, the sequence of ergodic averages
\begin{equation}\label{eq:ergaves}
\frac1n\sum_{k=0}^{n-1}\alpha^k(a)
\end{equation}
converges in norm as $n\to\infty$.

\item[(v)] We have
\begin{equation}
A^\alpha+\overline{\{a-\alpha(a)\mid a\in A\}}=A,
\end{equation}
where the closure is with respect to the norm topology.
\end{itemize}
These five statements imply the following statement:
\begin{itemize}
\item[(vi)] There exists a unique $\alpha$--invariant conditional expectation
$E$ from $A$ onto $A^\alpha$.
\end{itemize}
Furthermore, if (i)--(v) hold, then the conditional expectation $E$ in~(vi) is given
by the formula
\begin{equation}\label{eq:E}
E(a)=\lim_{n\to\infty}\frac1n\sum_{k=0}^{n-1}\alpha^k(a).
\end{equation}
\end{thm}

\begin{defi}
We say $\alpha$ is {\em uniquely ergodic relative to its fixed--point subalgebra}
if the equivalent statements (i)--(v) hold.
\end{defi}

\begin{proof}[Proof of Theorem~\ref{thm:ue}]
(i)$\implies$(ii) is clear.

\medskip\noindent
(ii)$\implies$(iii):
Suppose, to obtain a contradiction, that (ii) holds but
$x\in A$ and
\begin{equation}
x\notin\overline{A^\alpha+\{a-\alpha(a)\mid a\in A\}}.
\end{equation}
By the Hahn--Banach Theorem, there is a bounded linear functional $\phi:A\to\CC$ such that $\phi(x)\ne0$,
$\phi(A^\alpha)=\{0\}$ and $\phi\circ\alpha=\phi$.
Taking the real and imaginary parts, we may without loss of generality assume that $\phi$ is self--adjoint.
Let $\phi=\phi_+-\phi_-$ be the Jordan decomposition of $\phi$.
Then $\phi_+$ and $\phi_-$ are $\alpha$--invariant, by Observation~\ref{obs:Jordan}.
Moreover, $\phi_+$ and $\phi_-$ agree on $A^\alpha$.
Either both restrict to zero on $A^\alpha$, in which case $\phi_\pm(1)=0$ and $\phi_\pm=0$,
or $\phi_\pm$ are nonzero multiples of states on $A$
and by statement~(ii), $\phi_+$ and $\phi_-$ must agree on all of $A$.
This contradicts $\phi(x)\ne0$.

\medskip\noindent
(iii)$\implies$(iv):
Let $a\in A$ and $\eps>0$.
Let $c\in A^\alpha$ and $b\in A$ be such that
\begin{equation}
\|a-(c+b-\alpha(b))\|<\eps.
\end{equation}
If $n\ge m$, then
\begin{align}
\bigg\|\frac1n&\sum_{k=0}^{n-1}\alpha^k(a)-\frac1m\sum_{k=0}^{m-1}\alpha^k(a)\bigg\| \\
&<2\eps+\bigg\|\frac1n\sum_{k=0}^{n-1}\alpha^k(b-\alpha(b))
 -\frac1m\sum_{k=0}^{m-1}\alpha^k(b-\alpha(b))\bigg\| \\
&=2\eps+\bigg\|\frac1n(b-\alpha^n(b))+\frac1m(b-\alpha^m(b))\bigg\| \\
&\le2\eps+\frac{4\|b\|}m. \label{eq:3eps}
\end{align}
Taking $m$ large enough, the upper bound~\eqref{eq:3eps} can be made $<3\eps$.
This shows that the sequence of ergodic averages~\eqref{eq:ergaves} is Cauchy.

\medskip\noindent
(iv)$\implies$(vi)$+$\eqref{eq:E}:
Let $E$ be defined by the formula~\eqref{eq:E}.
Clearly, $E$ restricts to the identity map on $A^\alpha$.
One easily shows $\|E\|=1$ and $E\circ\alpha=\alpha\circ E=E$.
So $E$ is an $\alpha$--invariant conditional expectation from $A$ onto $A^\alpha$.
If $E':A\to A^\alpha$ is any $\alpha$--invariant conditional expectation onto $A^\alpha$,
then
\begin{equation}
E'(a)=\frac1n\sum_{k=0}^{n-1}E'(\alpha^k(a))=E'\bigg(\frac1n\sum_{k=0}^{n-1}\alpha^k(a)\bigg).
\end{equation}
Taking the limit as $n\to\infty$ gives
\begin{equation}
E'(a)=E'(E(a))=E(a).
\end{equation}

\medskip\noindent
(iv)$+$(vi)$+$\eqref{eq:E}$\implies$(i):
Let $\tau:A^\alpha\to\CC$ be a bounded linear functional.
Then $\tau\circ E$ is an $\alpha$--invariant extension of $\tau$ to $A$.
To show uniqueness, suppose $\phi:A\to\CC$ is any bounded, $\alpha$--invariant, linear
extension of $\tau$.
Then
\begin{equation}
\phi(a)=\frac1n\sum_{k=0}^{n-1}\phi(\alpha^k(a))=\phi\bigg(\frac1n\sum_{k=0}^{n-1}\alpha^k(a)\bigg).
\end{equation}
Taking the limit as $n\to\infty$ gives
\begin{equation}
\phi(a)=\phi(E(a))=\tau(E(a)),
\end{equation}
so $\phi=\tau\circ E$.

\medskip\noindent
We have now proved the equivalence of (i)--(iv), and that these imply (vi) and~\eqref{eq:E}.

\medskip\noindent
(i)$+$(vi)$\implies$(v):
Since $A=A^\alpha+\ker E$, it will suffice to show
\begin{equation}
\ker E\subseteq\overline{\{a-\alpha(a)\mid a\in A\}}.
\end{equation}
Suppose, to obtain a contradiction, $x\in\ker E$ but $x\notin\overline{\{a-\alpha(a)\mid a\in A\}}$.
By the Hahn--Banach Theorem, there is a bounded linear functional $\phi:A\to\CC$ such that
$\phi(x)\ne0$ and $\phi\circ\alpha=\phi$.
By~(i), we must have $\phi=\phi\circ E$, so $\phi(x)=0$, a contradiction.

\medskip\noindent
(v)$\implies$(iii) is clear.
\end{proof}

\begin{ques}\label{qn}
In Theorem~\ref{thm:ue}, is (vi) equivalent to (i)--(v)?
\end{ques}

It was kindly pointed out to us by Thierry Fack that the argument used in section \ref{sec:CFinf} applies more generally.
Indeed, as the following proposition shows,
the argument carries over to groups with property (RD), as defined by Jolissaint in \cite{J}.
Note that by \cite{dlH} this  includes the case of Gromov's hyperbolic groups.

\begin{prop}\label{prop:Fack}
Let $G$ be a group with property (RD) for a length function $L$
and let $\beta$ be an $L$--preserving automorphism of $G$ such that
all orbits of $\beta$ are either singletons or infinite.
Let $H=\{h\in G\mid\beta(h)=h\}$.
Then the automorphism $\alpha$ induced by $\beta$ on $C^*_r(G)$ is uniquely ergodic
relative to its fixed--point subalgebra,
which is the canonical copy of $C^*_r(H)$ in $C^*_r(G)$.  
\end{prop}
\begin{proof}
If $g\in G$ is such that $\beta(g)\neq g$, then
by \cite[Remark 1.2.2]{J} there exist positive numbers $C$ and $s$ such that
\begin{equation}
\bigg\|\frac{1}{n}\sum_0^{n-1}\alpha^k(\lambda_g)\bigg\|
\le C \bigg\|\frac{1}{n}\sum_0^{n-1}\alpha^k(\lambda_g)\bigg\|_{2,s,L}=\frac{C}{\sqrt n}(1+L(g))^s,
\end{equation}
and this upper bound approaches zero as $n$ goes to $\infty$.
If $\beta(g)=g$, then 
\begin{equation}
\frac{1}{n}\sum_0^{n-1}\alpha^k(\lambda_g)=\lambda_g
\end{equation}
for all $n$.
Now one easily sees that condition~(iv) of Theorem~\ref{thm:ue} holds and 
that $C^*_r(H)$ is the fixed--point subalgebra for $\alpha$.
\end{proof}

\section{The construction of reduced amalgamated free product C$^*$--algebras}
\label{sec:constr}

In this section we will review in some detail
and thereby set some notation for
the reduced amalgamated free product of C$^*$--algebras, which
was invented by Voiculescu~\cite{voic}.

We first set some notation concerning a right Hilbert \cstar--module $E$
over a \cstar--algebra $B$ (see \cite{la} for a general reference on Hilbert \cstar--modules).
If $x\in E$, then we let
\begin{equation}
|x|=\langle x,x\rangle^{1/2}\in B
\end{equation}
and the norm of $x$ is defined by 
\begin{equation}
\|x\|_E=\|\,|x|\,\|_B.
\end{equation}

Let $B$ be a unital \cstar--algebra,
let $I$ be a set with at least two elements
and for every $i\in I$ let $A_i$ be a unital \cstar--algebra
containing a copy of $B$ as a unital \cstar--subalgebra and having a conditional expectation 
$\phi_i:A_i\rightarrow B$ such that for each $a_i\in A_i$ there exists $x\in A_i$
for which $\phi_i(x^*a_i^*a_ix)\neq 0$.
We denote by $E_i=L^2(A_i,\phi_i)$ the right Hilbert \cstar--module over $B$ obtained
by separation and completion of $A_i$ with respect to the inner product
$\langle x,y\rangle=\phi_i(x^*y)$.
For $a_i\in A_i$,
we denote by $\ha_i$ the image of $a_i$ in $E_i$ under the canonical map.
There is a faithful $*$--representation $\pi_i$ of $A_i$ on $E_i$ by adjointable operators given by
\begin{equation}
\pi_i(x)(\hat{y})=(xy)\hat{\;},
\end{equation}
for $x,y\in A_i$.
We will often omit the reference to $\pi_i$
and write simply $av$ to denote $\pi_i(a)(v)$, for $a\in A_i$ and $v\in E_i$. 

This inclusion $B\subseteq A_i$ yields a copy of $B$ as a complemented Hilbert \cstar--submodule
of $E_i$,
and we write $E_i=B\oplus E_i^\circ$
and let $H_i:E_i\to E_i^\circ$ be the orthogonal projection onto $E_i^\circ$.
So, for example, we have
\begin{equation}
H_i(\ha)=(a-\phi_i(a))\hat{\;},\quad(a\in A_i).
\end{equation}
Since $\pi_i(b)$ sends $E_i^\circ$ into $E_i^\circ$ whenever $b\in B$, we
regard $E_i^\circ$ as equipped with a left $B$--action via $\pi_i$.
We consider the right Hilbert $B$--module
\begin{equation}\label{eq:Esum}
E= B\oplus \bigoplus_{\substack{m\in \NN \\i_1,\ldots,i_m\in I\\ i_j\neq i_{j+1}}}
  E^\circ_{i_1}\otimes_B E^\circ_{i_2}\otimes_B\cdots\otimes_B E^\circ_{i_m},
\end{equation}
where the tensor products are with respect to the right Hilbert $B$--module structures and the left actions
of $B$ described above,
and where the summand $B$ in~\eqref{eq:Esum}
denotes the \cstar--algebra $B$ with its usual Hilbert \cstar--module structure over itself.
There is a faithful $*$--representation of $A_i$ by adjointable operators on $E$, which
is denoted by $a\mapsto\lambda_a^i$
and which can be defined by
\begin{equation}\label{eq:lambdaaib}
\lambda_a^i(b)=\phi_i(ab)+H_i((ab)\hat{\;})\in B\oplus E^\circ_{i}, \quad(b\in B)
\end{equation}
and, considering a simple tensor
\begin{equation}\label{eq:xtdtx}
x_1\tdt x_m
\end{equation}
where $m\ge1$, $x_j\in E^\circ_{i_j}$, $i_1,\ldots,i_m\in I$ and $i_j\neq i_{j+1}$ for all $i=1,\ldots,m-1$,
by
\begin{equation}\label{eq:lambdaaixtx}
\lambda_a^i(x_1\tdt x_m)=
\begin{cases}
\begin{aligned}[b]
&H_i(\ha)\otimes x_1\tdt x_m \\
& +\phi_i(a)x_1\otimes x_2\tdt x_m,
\end{aligned}
&i\neq i_1\\[2ex]
\begin{aligned}[b]
&H_i(ax_1)\otimes x_2\tdt x_m \\
& +\langle (a^*)\hat{\;},x_1\rangle x_2\tdt x_m,
\end{aligned}
&i=i_1.
\end{cases}
\end{equation}
Note that for $b\in B$, $\lambda_b^i$ is the same for all $i$.
We will write $\lambda_a$ or simply $a$ instead of $\lambda^i_a$, when no confusion will result.

The reduced amalgamated free product \cstar--algebra
\begin{equation}
(A,\phi)={(*_B)}_{i\in I}(A_i,\phi_i)
\end{equation}
consists of the \cstar--algebra $A$ generated in ${\mathcal L}(E)$ by
the set $\{\lambda_a^i:a\in A_i,\ i\in I\}$ and the conditional expectation
$\phi:A\rightarrow B$ defined by
\begin{equation}
\phi(a)=\langle a 1_B,1_B\rangle,\quad(a\in A).
\end{equation}
Thus, the \cstar--algebra $A$ is the closed span of $B$
together with the set of all words of the form
\begin{equation}\label{eq:w}
w=a_1\dots a_n
\end{equation}
where $a_i\in A^\circ_{k(i)}$, $k(1),\ldots,k(n)\in I$ and $k(i)\neq k(i+1)$
for all $i\in\{1,\dots,n-1\}$.

\section{Some norm estimates in reduced amalgamated free product C$^*$--algebras}
\label{sec:normest}

The main result of this section is the following norm estimate,
which applies to certain linear combinations
of words of length $n$ in reduced amalgamated free product \cstar--algebras.
It is a version of the Haagerup inequality.

\begin{prop}\label{prop:ineq}
Suppose $n\ge1$ and consider
\begin{equation}
f=\sum_{k\in\fk}a_{k,1}a_{k,2}\cdots a_{k,n}\in A,
\end{equation}
where $\fk$ is a finite subset of $I^n$ such that for all $k=(k(1),\ldots,k(n))\in\fk$ we have
$k(i)\ne k(i+1)$ for all $i\in\{1,\ldots,n-1\}$ and where $a_{k,i}\in A_{k(i)}^\circ$ for all
$k\in\fk$ and $i\in\{1,\ldots,n\}$.
Suppose, furthermore, that
\begin{equation}\label{eq:knek'}
\text{if }k,k'\in\fk\text{ and }\,k\ne k',\text{ then }k(1)\ne k'(1)\text{ and }k(n)\ne k'(n).
\end{equation}
Then
\begin{equation}\label{eq:fnm}
\|f\|\le(2n+1)\bigg(\sum_{k\in\fk}\prod_{i=1}^n\|a_{k,i}\|^2\bigg)^{1/2}.
\end{equation}
\end{prop}

Before we get to the proof, we consider some preliminary constructions and results.
Let us define some elementary adjointable operators on $E$, in terms of which we
will later describe the action of a word $w$ as in~\eqref{eq:w} on a tensor $x_1\tdt x_m$
in~\eqref{eq:xtdtx}.

\begin{notation}
Let $P_0$ denote the orthogonal projection of $E$ onto the summand $B\subseteq E$ and for $m\ge1$
let $P_m$ denote the orthogonal projection of $E$ onto
\begin{equation}
\bigoplus_{\substack{i_1,\ldots,i_m\in I\\ i_j\neq i_{j+1}}}
  E^\circ_{i_1}\otimes_B E^\circ_{i_2}\otimes_B\cdots\otimes_B E^\circ_{i_m}.
\end{equation}
\end{notation}

\begin{notation}
For $k\in I$, let $Q_k$ denote the orthogonal projection of $E$ onto
\begin{equation}
\bigoplus_{\substack{m\ge1 \\i_1,\ldots,i_m\in I\\ i_j\neq i_{j+1} \\ i_1=k}}
  E^\circ_{i_1}\otimes_B E^\circ_{i_2}\otimes_B\cdots\otimes_B E^\circ_{i_m}.
\end{equation}
Note that $Q_k$ and $P_m$ commute.
\end{notation}

\begin{notation}
Given $k\in I$ and $y\in E_k^\circ$, let $\psi(y)=\psi_k(y)\in\fl(E)$ be given by
\begin{equation}\label{eq:psiyb}
\psi(y)b=(yb)\hat{\;}\in E_k^\circ,\quad(b\in B)
\end{equation}
and, for $x_1\tdt x_m$ as in~\eqref{eq:xtdtx},
\begin{equation}\label{eq:psiyxtx}
\psi(y)(x_1\tdt x_m)=
\begin{cases}
0,&i_1=k \\
y\otimes x_1\tdt x_m,& i_1\ne k.
\end{cases}
\end{equation}
Therefore, we have
\begin{gather}
\psi(y)=Q_k\psi(y)(1-Q_k), \displaybreak[2] \\[1ex]
\psi(y)^*b=0,\quad (b\in B) \label{eq:psiy*b} \displaybreak[2] \\[1ex]
\psi(y)^*(x_1\tdt x_m)=
\begin{cases}
0,&i_1\ne k \\
\langle y,x_1\rangle,&i_1=k,\,m=1 \\
\langle y,x_1\rangle x_2\otimes x_3\tdt x_m,& i_1=k,\,m>1
\end{cases} \label{eq:psiy*xtx} \displaybreak[2] \\[1ex]
\psi(y)^*\psi(y)=|y|^2(1-Q_k) \displaybreak[2] \\[1ex]
\|\psi(y)\|=\|y\|. \label{eq:psiynorm}
\end{gather}
\end{notation}

\begin{notation}
For $k\in I$ and $a\in A_k$, we let $\rho(a)=\rho_k(a)\in\fl(E)$ be defined by
\begin{equation}\label{eq:rhoab}
\rho(a)b=0,\quad(b\in B)
\end{equation}
and, for $x_1\tdt x_m$ as in~\eqref{eq:xtdtx},
\begin{equation}\label{eq:rhoaxtx}
\rho(a)(x_1\tdt x_m)=
\begin{cases}
(H_k(ax_1))\otimes x_2\tdt x_m,&i_1=k \\
0,&i_1\ne k.
\end{cases}
\end{equation}
(Recall that $H_k:E_k\to E_k^\circ$ is the orthogonal projection.)
Therefore, we have
\begin{gather}
\rho(a)=Q_k\rho(a)Q_k \\[1ex]
\|\rho(a)\|\le\|a\|. \label{eq:rhoanorm}
\end{gather}
\end{notation}

To ease notation, for $a\in A_k$ we let
\begin{equation}
\ha^\dag=(a^*)\hat{\;}\in E_k.
\end{equation}

The following lemma describes how a word $w=a_1\cdots a_n$ as in~\eqref{eq:w} can act
on a tensor $x_1\tdt x_m$ as in~\eqref{eq:xtdtx}.
What can happen is: $w$ can first devour some initial string $x_1\tdt x_q$ of the tensor.
Then it can either push some more stuff onto the tensor from the left, or it can instead
act on the next letter $x_{q+1}$ and then push some more stuff onto the tensor from the left.
This is all that can happen, because neighboring letters in $w$ and neighboring $x_j$ in $x_1\tdt x_m$
are constrained to come from different $A_k^\circ$, respectively different $E_i^\circ$.
It's not too difficult to see this by considering some examples.
We give a more precise statement and a rigorous proof below.

\begin{lemma}\label{lem:PwP}
Let $n\ge1$ and let $k=(k(1),\ldots,k(n))\in I^n$ be such that $k(i)\ne k(i+1)$ for all $i\in\{1,\ldots,n-1\}$.
Let $w=a_1\cdots a_n$, where $a_i\in A_{k(i)}^\circ$ for all $i\in\{1,\ldots,n\}$.
Let $m,r\ge0$ be integers.
\begin{itemize}
\item[(i)] If $r>m+n$ or $r<|m-n|$, then $P_rwP_m=0$.
\item[(ii)] If $r=m+n-2s$ with $s\in\{0,1,\ldots,\min(m,n)\}$,
then
\begin{equation}
\begin{aligned}
P_rwP_m=\psi(\ha_1)&\psi(\ha_2)\cdots\psi(\ha_{n-s}) \\
 &\cdot\psi(\ha_{n-s+1}^\dag)^*\psi(\ha_{n-s+2}^\dag)^*\cdots\psi(\ha_n^\dag)^*P_m.
\end{aligned}
\end{equation}
\item[(iii)] If $r=m+n-2s+1$ with $s\in\{1,2,\ldots,\min(m,n)\}$, then
\begin{equation}
\begin{aligned}
P_rwP_m&=\psi(\ha_1)\psi(\ha_2)\cdots\psi(\ha_{n-s}) \\
&\quad\quad\cdot\rho(a_{n-s+1})\,
 \psi(\ha_{n-s+2}^\dag)^*\psi(\ha_{n-s+3}^\dag)^*\cdots\psi(\ha_n^\dag)^*P_m.
\end{aligned}
\end{equation}
\end{itemize}
\end{lemma}
\begin{proof}
The following equation is equivalent to parts~(i)--(iii) of the Lemma taken together:
\begin{equation}\label{eq:bigsum}
\begin{aligned}
wP_m&=\sum_{s=0}^{\min(m,n)}
\begin{aligned}[t]
P_{n+m-2s}&\psi(\ha_1)\cdots\psi(\ha_{n-s}) \\
 &\cdot\psi(\ha_{n-s+1}^\dag)^*\cdots\psi(\ha_n^\dag)^*P_m
\end{aligned} \\
&\;+\sum_{s=1}^{\min(m,n)}
\begin{aligned}[t]
P_{n+m-2s+1}&\psi(\ha_1)\cdots\psi(\ha_{n-s})\rho(a_{n-s+1}) \\
 &\cdot\psi(\ha_{n-s+2}^\dag)^*\cdots\psi(\ha_n^\dag)^*P_m.
\end{aligned}
\end{aligned}
\end{equation}
We will prove~\eqref{eq:bigsum} by induction on $n$.
For $n=1$, taking first $m\ge1$ and using the fact that $\phi_{k(1)}(a_1)=0$
together with~\eqref{eq:lambdaaixtx}, \eqref{eq:psiyxtx},
\eqref{eq:psiy*xtx}, and~\eqref{eq:rhoaxtx}, we find
\begin{align}
a_1P_m&=\big(\psi(\ha_1)+\rho(a_1)+\psi(\ha_1^\dag)^*\big)P_m \\
&=P_{m+1}\psi(\ha_1)P_m+P_m\rho(a_1)P_m+P_{m-1}\psi(\ha_1^\dag)^*P_m, \label{eq:Pa1P}
\end{align}
while in the case $m=0$, using~\eqref{eq:lambdaaib}, \eqref{eq:psiyb},
\eqref{eq:psiy*b}, and~\eqref{eq:rhoab}, we find
\begin{equation}\label{eq:a1P0}
a_1P_0=\psi(\ha_1)P_0=P_1\psi(\ha_1)P_0.
\end{equation}
Thus,~\eqref{eq:bigsum} is proved in the case $n=1$.

Now let $n\ge2$ and set $w'=a_2a_3\cdots a_n$.
By the induction hypothesis, we have
\begin{align}
w'P_m&=\sum_{s=0}^{\min(m,n-1)}
\begin{aligned}[t]
P_{n+m-2s-1}&\psi(\ha_2)\cdots\psi(\ha_{n-s}) \\
 &\cdot\psi(\ha_{n-s+1}^\dag)^*\cdots\psi(\ha_n^\dag)^*P_m
\end{aligned} \label{eq:w'sum1} \\
&\;+\sum_{s=1}^{\min(m,n-1)}
\begin{aligned}[t]
P_{n+m-2s}&\psi(\ha_2)\cdots\psi(\ha_{n-s})\rho(a_{n-s+1}) \\
 &\cdot\psi(\ha_{n-s+2}^\dag)^*\cdots\psi(\ha_n^\dag)^*P_m.
\end{aligned} \label{eq:w'sum2} 
\end{align}
Now we multiply both sides of~\eqref{eq:w'sum1} and~\eqref{eq:w'sum2}
on the left by $a_1$, and use~\eqref{eq:Pa1P}
and~\eqref{eq:a1P0}, as needed.
For example, from~\eqref{eq:w'sum1} consider
\begin{equation}\label{eq:a1Pstuff}
a_1P_{n+m-2s-1}\psi(\ha_2)\cdots\psi(\ha_{n-s}) \\
 \psi(\ha_{n-s+1}^\dag)^*\cdots\psi(\ha_n^\dag)^*P_m.
\end{equation}
If $s<n-1$, then the initial part of~\eqref{eq:a1Pstuff} is
\begin{equation}
\begin{aligned}
a_1P_{n+m-2s-1}\psi(\ha_2)&=P_{n+m-2s}\psi(\ha_1)P_{n+m-2s-1}\psi(\ha_2) \\
&\qquad+P_{n+m-2s-1}\rho(a_1)P_{n+m-2s-1}\psi(\ha_2) \\
&\qquad+P_{n+m-2s-2}\psi(\ha_1^\dag)^*P_{n+m-2s-1}\psi(\ha_2) \\
&=P_{n+m-2s}\psi(\ha_1)P_{n+m-2s-1}\psi(\ha_2) \\
&=P_{n+m-2s}\psi(\ha_1)\psi(\ha_2),
\end{aligned}
\end{equation}
where we have used
\begin{gather}
\rho(a_1)P_{n+m-2s-1}\psi(\ha_2)
 =\rho(a_1)Q_{k(1)}P_{n+m-2s-1}Q_{k(2)}\psi(\ha_2)=0 \\
\psi(\ha_1^\dag)^*P_{n+m-2s-1}\psi(\ha_2)
 =\psi(\ha_1^\dag)^*Q_{k(1)}P_{n+m-2s-1}Q_{k(2)}\psi(\ha_2)=0.
\end{gather}
If $s=n-1<m$, then the initial part of~\eqref{eq:a1Pstuff} is
\begin{equation}
\begin{aligned}
a_1P_{m-s}\psi(\ha_2^\dag)^*&=P_{m-s+1}\psi(\ha_1)P_{m-s}\psi(\ha_2^\dag)^* \\
&\qquad+P_{m-s}\rho(a_1)P_{m-s}\psi(\ha_2^\dag)^* \\
&\qquad+P_{m-s-1}\psi(\ha_1^\dag)^*P_{m-s}\psi(\ha_2^\dag)^* \\
&=P_{m-s+1}\psi(\ha_1)\psi(\ha_2^\dag)^*+P_{m-s}\rho(a_1)\psi(\ha_2^\dag)^* \\
&\qquad+P_{m-s-1}\psi(\ha_1^\dag)^*\psi(\ha_2^\dag)^*,
\end{aligned}
\end{equation}
while if $s=n-1=m$, then the initial part of~\eqref{eq:a1Pstuff} is
\begin{equation}
a_1P_0\psi(\ha_2^\dag)^*=P_1\psi(\ha_1)P_0\psi(\ha_2^\dag)^*=P_1\psi(\ha_1)\psi(\ha_2^\dag)^*.
\end{equation}
Turning now to~\eqref{eq:w'sum2}, we consider
\begin{equation}\label{eq:a1Pstuff2}
a_1P_{n+m-2s}\psi(\ha_2)\cdots\psi(\ha_{n-s})\rho(a_{n-s+1})
 \psi(\ha_{n-s+2}^\dag)^*\cdots\psi(\ha_n^\dag)^*P_m.
\end{equation}
We find that the initial part of~\eqref{eq:a1Pstuff2} is
\begin{equation}
a_1P_{n+m-2s}\psi(\ha_2)=
\begin{cases}
P_{n+m-2s+1}\psi(\ha_1)\psi(\ha_2),&s<n-1 \\
P_{m-s+2}\psi(\ha_1)\rho(a_2),&s=n-1.
\end{cases}
\end{equation}
Taking all of these cases into account, we prove~\eqref{eq:bigsum}.
\end{proof}

\begin{lemma}\label{lem:PrfPm}
Let $f$ be as in Proposition~\ref{prop:ineq}.
Let $m,r$ be nonnegative integers.
Then
\begin{equation}\label{eq:PfPnm}
\|P_rfP_m\|^2\le\sum_{k\in\fk}\prod_{i=1}^n\|a_{k,i}\|^2.
\end{equation}
\end{lemma}
\begin{proof}
If $r<|m-n|$ or $r>m+n$, then by Lemma~\ref{lem:PwP}(i), we have $P_rfP_m=0$.

\noindent{\em Case I.}
Suppose $r=m+n-2s$ for $s\in\{0,1,\ldots,\min(m,n)\}$ and with $s<n$.
By Lemma~\ref{lem:PwP}(ii), we have
\begin{equation}\label{eq:PrfPm}
P_rfP_m=\sum_{k\in\fk}\psi(\ha_{k,1})\cdots\psi(\ha_{k,n-s})
 \psi(\ha_{k,n-s+1}^\dag)^*\cdots\psi(\ha_{k,n}^\dag)^*P_m
\end{equation}
and
\begin{multline}
(P_rfP_m)^*(P_rfP_m)= \\
=\sum_{k,k'\in\fk}
\begin{aligned}[t]
P_m
&\psi(\ha_{k,n})\cdots\psi(\ha_{k,n-s+1})
 \psi(\ha_{k,n-s}^\dag)^*\cdots\psi(\ha_{k,1}^\dag)^* \\
&\cdot\psi(\ha_{k',1})\cdots\psi(\ha_{k',n-s})
 \psi(\ha_{k',n-s+1}^\dag)^*\cdots\psi(\ha_{k',n}^\dag)^*P_m.
\end{aligned}
\end{multline}
By the hypothesis~\eqref{eq:knek'}, if $k\ne k'$, then $k(1)\ne k'(1)$ and, consequently,
\begin{equation}
\psi(\ha_{k,1}^\dag)^*\psi(\ha_{k',1})=\psi(\ha_{k,1}^\dag)^*Q_{k(1)}Q_{k'(1)}\psi(\ha_{k',1})=0.
\end{equation}
Therefore, using also~\eqref{eq:psiynorm}, we get
\begin{equation}\label{eq:PfPnm1}
\|P_rfP_m\|^2\le
\sum_{k\in\fk}\prod_{i=1}^n\|\ha_{k,i}\|^2
\le\sum_{k\in\fk}\prod_{i=1}^n\|a_{k,i}\|^2.
\end{equation}

\noindent{\em Case II.}
Suppose $r=m+n-2s$ for $s=n\le m$.
Then~\eqref{eq:PrfPm} becomes
\begin{equation}
P_rfP_m=\sum_{k\in\fk}\psi(\ha_{k,1}^\dag)^*\cdots\psi(\ha_{k,n}^\dag)^*P_m
\end{equation}
and we have
\begin{equation}
(P_rfP_m)(P_rfP_m)^*=\sum_{k,k'\in\fk}
\psi(\ha_{k,1}^\dag)^*\cdots\psi(\ha_{k,n}^\dag)^*P_m\psi(\ha_{k',n})\cdots\psi(\ha_{k',1}).
\end{equation}
Again, by the hypothesis~\eqref{eq:knek'}, if $k\ne k'$, then $k(n)\ne k'(n)$ and, consequently,
\begin{equation}
\psi(\ha_{k,n}^\dag)^*P_m\psi(\ha_{k',n})=\psi(\ha_{k,n}^\dag)^*Q_{k(1)}P_mQ_{k'(1)}\psi(\ha_{k',n})=0.
\end{equation}
Using again~\eqref{eq:psiynorm}, we get~\eqref{eq:PfPnm1} also in this case.

\vskip1ex
\noindent{\em Case III.}
Suppose $r=m+n-2s+1$ for $s\in\{1,\ldots,\min(m,n)\}$.
Then using Lemma~\ref{lem:PwP}(iii) and proceeding
similarly to above, we obtain the estimate
\begin{equation}\label{eq:PfPnm2}
\|P_rfP_m\|^2\le\sum_{k\in\fk}\|a_{k,n-s+1}\|^2\prod_{\substack{1\le i\le n \\ i\ne n-s+1}}\|\ha_{k,i}\|^2
\le\sum_{k\in\fk}\prod_{i=1}^n\|a_{k,i}\|^2.
\end{equation}
\end{proof}

\begin{rk}
The left--hand inequalities in~\eqref{eq:PfPnm1} and~\eqref{eq:PfPnm2} are
better than required in~\eqref{eq:PfPnm}.
In fact,~\eqref{eq:PfPnm1} and~\eqref{eq:PfPnm2} seem to be quite close in spirit to
the inequality obtained in~\cite[3.3]{DHR}, which applied to free products with amalgamation
over the scalars.
\end{rk}

\begin{proof}[Proof of Proposition~\ref{prop:ineq}]
Let $\sigma:B\to\fl(\HEu)$ be a faithful $*$--representation of $B$ on a Hilbert space $\HEu$.
Then the internal tensor product $\HEut=E\otimes_\sigma\HEu$ is a Hilbert space
and the $*$--representation $\sigmat:\fl(E)\to\fl(\HEut)$ given by $\sigmat(a)=a\otimes\id_\HEu$
is faithful.

Let $v\in\HEut$.
Then
\begin{equation}\label{eq:Prfv}
\|\sigmat(f)v\|^2=\sum_{r=0}^\infty\|\sigmat(P_rf)v\|^2.
\end{equation}
Let
\begin{equation}
\gamma=\bigg(\sum_{k\in\fk}\prod_{i=1}^n\|a_{k,i}\|^2\bigg)^{1/2}.
\end{equation}
Then
\begin{align}
\|\sigmat(P_rf)v\|^2&=\bigg\|\sum_{m=|r-n|}^{r+n}\sigmat(P_rfP_m)v\bigg\|^2 \label{eq:Prf1} \\
&\le\bigg(\sum_{m=|r-n|}^{r+n}\|\sigmat(P_rfP_m)v\|\bigg)^2 \displaybreak[2] \\
&\le\bigg(\sum_{m=|r-n|}^{r+n}\gamma\|\sigmat(P_m)v\|\bigg)^2 \label{eq:Prf3} \displaybreak[2] \\
&\le\bigg(\sum_{m=|r-n|}^{r+n}\gamma^2\bigg)\bigg(\sum_{m=|r-n|}^{r+n}\|\sigmat(P_m)v\|^2\bigg)
 \label{eq:Prf4} \displaybreak[2] \\
&\le(2n+1)\gamma^2\sum_{m=|r-n|}^{r+n}\|\sigmat(P_m)v\|^2 \label{eq:Prf5}
\end{align}
where we used Lemma~\ref{lem:PwP}(i) to get~\eqref{eq:Prf1},
Lemma~\ref{lem:PrfPm} to get~\eqref{eq:Prf3}
and the Cauchy--Schwarz inequality to get~\eqref{eq:Prf4}.
From~\eqref{eq:Prfv} and~\eqref{eq:Prf1}--\eqref{eq:Prf5}, we get
\begin{align}
\|\sigmat(f)v\|^2&\le(2n+1)\gamma^2\sum_{r=0}^\infty\sum_{m=|r-n|}^{r+n}\|\sigmat(P_m)v\|^2 \\
&=(2n+1)\sum_{m=0}^\infty\sum_{r=|m-n|}^{m+n}\|\sigmat(P_m)v\|^2 \\
&\le(2n+1)^2\gamma^2\sum_{m=0}^\infty\|\sigmat(P_m)v\|^2 \\
&=(2n+1)^2\gamma^2\|v\|^2.
\end{align}
This shows
$\|\sigmat(f)\|\le(2n+1)\gamma$, which implies~\eqref{eq:fnm}.
\end{proof}

\section{Free shifts}
\label{sec:freeshift}

Let $D$ be a unital \cstar--algebra,
and  let $E^D_B:D\rightarrow B$ be a conditional expectation onto a unital \cstar--subalgebra $B$.
For each integer $i\in\ZZ$ let $(A_i,\phi_i)$ be a copy of $(D,E^D_B)$.
Let 
\begin{equation}\label{eq:Aphi}
(A,\phi)={(*_B)}_{i\in I}(A_i,\phi_i)
\end{equation}
be the reduced amalgamated free product,
and let $a\mapsto\lambda_a^i$ denote the embedding of $A_i$ in the free product algebra $A$ arising from
the free product construction, as descibed in section~\ref{sec:constr}.
The free--shift automophism $\alpha$ on $A$ is the
automorphism of $A$
given by $\alpha(\lambda_a^i)=\lambda_a^{i+1}$ for all $a\in A$ and $i\in\ZZ$.

\begin{thm}\label{thm:afp}
Let $\alpha$ be the free--shift automorphism on the amalgamated free product C$^*$--algebra
$A$ as given in~\eqref{eq:Aphi} above.
Then $B$ is the fixed--point subalgebra for $\alpha$
and $\alpha$ is uniquely ergodic relative to its fixed--point subalgebra.
\end{thm}

\begin{proof}
We will show
\begin{equation}\label{eq:ergavephi}
\lim_{n\to\infty}\frac1n\sum_{k=0}^{n-1}\alpha^k(a)=\phi(a)
\end{equation}
for every $a\in A$.
This will imply that $B$ is the fixed--point subalgebra for $\alpha$ and that condition~(iv) of
Theorem~\ref{thm:ue} holds.

It will suffice to show~\eqref{eq:ergavephi} for all $a\in B$ and words $a$ of the form
$w=a_1a_2\cdots a_p$ for some $p\ge1$ and $a_i\in A_{k(i)}^\circ$,
and some $k(i)\in\ZZ$ with $k(i)\ne k(i+1)$, $i=1,\ldots,p-1$.
Since $B$ is $\alpha$ invariant,~\eqref{eq:ergavephi} is clear for $a\in B$.
So assume $a=w$ as above.
Then $\phi(w)=0$ and
$\sum_{k=0}^{n-1}\alpha^k(w)$ is a finite linear combination of words of length $p$
to which Proposition~\ref{prop:ineq} applies, and we have
\begin{equation}
\bigg\|\frac1n\sum_{k=0}^{n-1}\alpha^k(w)\bigg\|\le\frac1n(2p+1)n^{1/2}\prod_{i=1}^p\|a_i\|.
\end{equation}
Thus, we get
\begin{equation}
\lim_{n\to\infty}\bigg\|\frac1n\sum_{k=0}^{n-1}\alpha^k(w)\bigg\|=0,
\end{equation}
as required.
\end{proof}

\end{document}